\documentclass{svproc}
\usepackage{url}

\usepackage{amssymb,latexsym,amsmath,amsfonts} 
\usepackage{graphicx,psfrag,overpic,contour,color} 
\usepackage{enumitem}

\newcommand\Lemref[1]{Lemma~\ref{#1}}

\newcommand\Figref[1]{Figure~\ref{#1}}

\def\ol#1{\overline{#1}}

\def\wkl{\,<\mskip-10mu)\mskip5mu}

\def\Sphere{\mathcal S^2}

\def\wh#1{\widehat{#1}}
\def\klzwi{\mskip 1mu}

\def\phi{\varphi}

\definecolor{blaucmyk}{cmyk}{1.00,0.40,0.00,0.20}

\definecolor{rotcmyk}{cmyk}{0.00,1.00,1.00,0.00}

\definecolor{hhgrau}{cmyk}{0.00,0.00,0.00,0.07}

\definecolor{hgrau}{cmyk}{0.00,0.00,0.00,0.15}

\definecolor{grau}{cmyk}{0.00,0.00,0.00,0.40}

\definecolor{dgrau}{cmyk}{0.00,0.00,0.00,0.85}

\definecolor{gruen}{cmyk}{1.00,0.00,0.90,0.30} 

\definecolor{gelb}{cmyk}{0.00,0.00,1.00,0.00} 

\definecolor{gelbb}{cmyk}{0.00,0.00,0.30,0.00} 

\definecolor{gelbbb}{cmyk}{0.00,0.00,0.15,0.00} 

\definecolor{weiss}{gray}{1}

\definecolor{lila}{rgb}{0.8,0.0,0.8}

\def\Square{\vbox{\hrule\hbox{\vrule height 1.5mm\hskip 1.5mm%
                         \vrule height 1.5mm}\hrule}}             
\newenvironment{Proof}{\noindent{\em Proof.}}
        {~\hfill\Square\par\smallskip}

\begin{document}
\mainmatter              
\title{Geometric analysis of Bennett's spherical\\ 8-bar linkage and its spatial counterpart}
\titlerunning{Geometric analysis of Bennett's spherical 8-bar linkage}  
\author{Hellmuth Stachel\inst{1}}
\authorrunning{Hellmuth Stachel} 
\institute{Vienna University of Technology, Vienna, Austria,\\
\email{stachel@dmg.tuwien.ac.at},\\ WWW home page:
\texttt{https://www.dmg.tuwien.ac.at/stachel/}
}

\maketitle              


\begin{abstract}
We provide a geometric approach to two combinatorically symmmetric overconstrained spatial linkages.
Both contain eight bodies and twelve revolute joints and collapse in aligned poses. 
The first one is spherical and the union of six spherical isograms.
It is the spherical image of a Bricard octahedron of type~3 and was already analysed  1912 by Bennett.
The second linkage is the dualized version and composed from six Bennett isograms. 
Our approach via line reflections discloses some symmetries at spatial poses. 
%
\end{abstract}
%
\section{Introduction}
R.\ Bricard's flexible octahedra \cite{Bricard} and G.~T.\ Bennett's isogram \cite{Bennett2} are fascinating examples of overconstrained spatial mechanisms with revolute joints.
Both play an important role at the two types of overconstrained linkages which are presented below.

\begin{figure}[htb] 
   \centerline{\includegraphics[width=105mm]{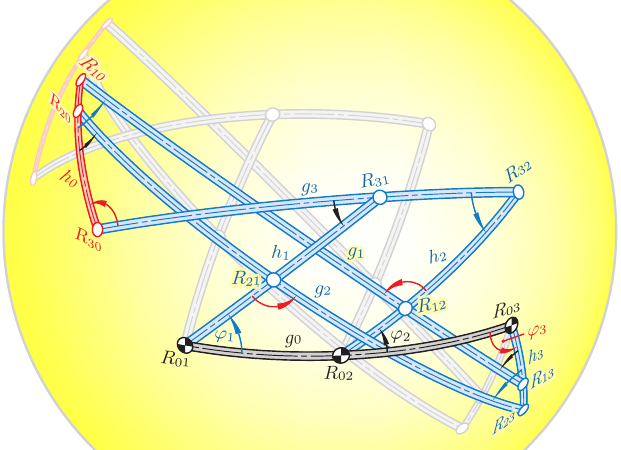}}
  \caption{Bennett's composition of spherical isograms presented in \cite[Figs.~5 and 6\,]{Bennett} is an overconstrained spherical linkage with 8 bars and 12 joints, arranged in six isograms.} 
  \label{fig:spher_linkage1}
\end{figure}

The spherical image of a Bricard octahedron of type~3 gives a spherical 8-bar linkage that is the union of six spherical isograms (\Figref{fig:spher_linkage1}).
Bennett proved in \cite{Bennett} with the help of trigonometric relations that this linkage is flexible and used this statement to confirm the flexibility of the type~3 octahedra.
A particular planar counterpart of this spherical linkage is attributed to G.\ Darboux \cite{Darboux}.

 E.\ Study's principle of transference offers a possibility to transfer statements from the geometry on the sphere to spatial geometry.
It is based on the dualization of quaternions (see, e.g., \cite{Sta_071}) and allows to transfer Bennett's spherical linkage to a spatial linkage with eight links and twelve hinges that is the union of six Bennett's isograms (\Figref{fig:spatial_linkage_kmplt}).

Also this linkage is not new.
E.\ Baker analysed in \cite[Figs.~8 and 11]{Baker1} the spherical linkage by providing closure equations in terms of trigonometric functions.\footnote{
It is interesting that in his paper Baker estimated Bennett's considerations in \cite{Bennett} as `inpenetrable' and `obscure'.}   
In the same paper, the planar and the spatial versions were discussed as well.
Moreover, Baker presented other compounds of Bennett isograms and revisited them along with flexible 6R-loops in a series of publications, e.g., in \cite{Baker2}.

\begin{figure}[htb] 
  \centerline{\includegraphics[width=48mm]{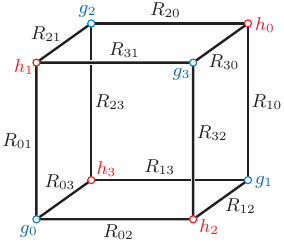}}
  \caption{At this linkgraph the links and joints of Bennet's 6-bar linkage (\Figref{fig:spher_linkage1}) and its spatial counterpart (\Figref{fig:spatial_linkage_kmplt}) are represented as knots and bars.}
  \label{fig:graph}
\end{figure}

The spatial version of Bennett's 8-bar linkage was also studied in a much more general context in \cite{Hege}.
Two different factorizations of a rational curve on the Study quadric, represented in terms of dual quaternions, give rise to closed revolute linkages.
In particular, the factorizations of rational cubics paves the way to the spatial 6-bar linkage (note the numeric Example~3 in \cite{Hege}) with the spherical version as the zero-offset type.
Their linkgraph, which represents the links and joints of the linkage as knots and bars of a graph (\Figref{fig:graph}), is equivalent to the 1-skeleton of a cube thus revealing the combinatorial symmetry of this structure.

\medskip
The goal of this article is to provide a novel geometric approach to Bennett's spherical 8-bar linkage and its spatial counterpart.
Thus, we are able to disclose some symmetries at their general poses.
In our approach, reflections in lines, i.e., rotations through $180^\circ$ (also called `halfturns') play a major role.
Several pictures are included in order to support the geometric reasoning and to illustrate the linkages.

%
\section{The spherical 8-bar linkage}
Let $\Sphere$ be the unit sphere with center $O$.
For better understanding, we begin with recalling some particularities of the geometry on the sphere $\Sphere$. 
Further details of spherical geometry can be found, e.g., in \cite[Sect.~10.1]{Conics}.

Since there are no straight lines on $\Sphere$, their role is taken over by {\em great circles}, i.e., circles in planes through the center $O$.  
If points of $\Sphere$ are located on the same great circle, we speak of {\em aligned} points. 

\begin{itemize}[itemsep=1pt]
\item Each point $P\in\Sphere$ has an {\em antipode} $P^\ast$, and each great circle passing through $P$ contains also $P^\ast$.
Two different great circles $g_1, g_2$ intersects at a pair of antipodal points.

\item For any two different points $P$ and $Q\ne P^\ast$ there exists a unique connecting great circle. 
When we speak of the {\em arc} or the {\em side} $PQ$ for any two non-antipodal points $P, Q$, we always mean on the connecting great circle the shorter of the two arcs terminated by $P$ and $Q$.

\item The {\em spherical distance} $\wh{PQ}$ is defined as length of the side $PQ$ measured on $\Sphere$.
It equals the central angle $\wkl POQ$ of the arc $PQ$ and obeys $0 \le \wh{PQ} \le \pi$, when including the limiting cases $Q=P$ and $Q=P^\ast$. 
We consider only spherical triangles and quadrangles with sides of lengths smaller than $\pi$.   


\item The {\em spherical bisector} of any two different points $P,Q\in \Sphere$ is defined as the set of points $X\in \Sphere$ obeying $\wh{XP} = \wh{XQ}$. 
Since equal spherical distances imply also equal 3D-distances along the chords, the bisector is the great circle in the plane of symmetry between $P$ and $Q$.

\item Given two great circles $g_1, g_2$, the {\em angle} $\wkl g_1g_2$ is defined as the angle between the tangent lines $t_1, t_2$ at their intersection points. 
 This angle equals the angle between the planes spanned by the circles and obeys $\wkl g_1g_2 = \wkl t_1t_2 \le \frac\pi 2$. 

\item The angle $\wkl QPR$ between the arcs $PQ$ and $PR$ on $\Sphere$ can be bounded by $0 \le \wkl QPR \le \pi$.
Angles of rotations about a point $P\in\Sphere$ can be signed by defining that  counter-clockwise rotations (when seen from outside) are positive. 
Note that the sign changes when $P$ is replaced by its antipode $P^\ast$.

\item The points of intersection between different great circles $g_1, g_2$ are the spherical centers of the unique {\em common perpendicular} great circle of $g_1$ and $g_2$.  
\end{itemize}

\begin{definition}\label{def:isogram}  \em
A spherical quadrangle $P_1\dots P_4$ is called a {\em spherical isogram} if opposite sides have the same lengths.
The quadrangle can be free of self-intersections or {\em crossed}.
In the latter case two opposite sides are intersecting. 
\end{definition}

\begin{figure}[hbt] 
  \centerline{\includegraphics[width=85mm]{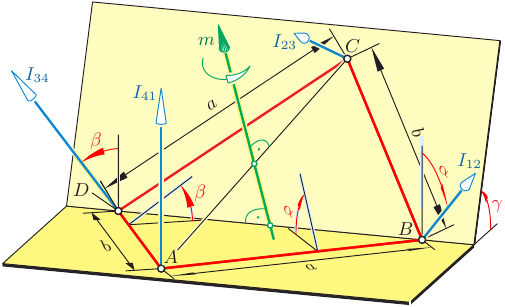}}
  \caption{At a skew isogram $ABCD$ with side lengths $a, b$ the dihedral angles $\alpha,\beta$ of its convex hull match the proportion $a:b = \sin\alpha : \sin\beta$.}
  \label{fig:Bennett_angles}
\end{figure}

\begin{figure}[hbt] 
  \centerline{\includegraphics[width=118mm]{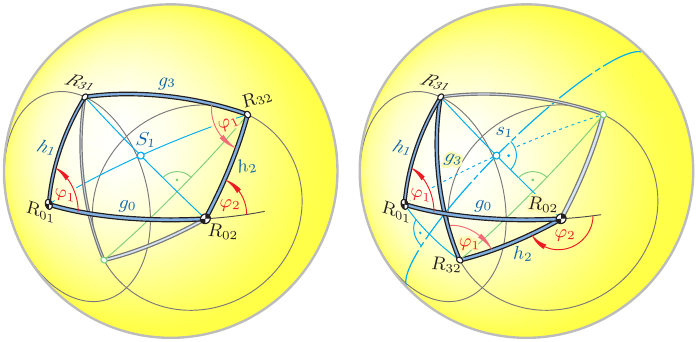}}
  \caption{At a non-crossed spherical isogram (left) the common midpoint $S_1$ of both diagonal arcs is a center of symmetry.
At a crossed spherical isogram (right) the axis $s_1$ of symmetry meets the two diagonal arcs orthogonally at there midpoints.
\\
If one side remains fixed, the coupler motion of a spherical isogram splits into two rational motions with bifurcations at the aligned poses.
} 
  \label{fig:sph_isogram}
\end{figure}

All skew isograms have an axis $m$ of symmetry.
It intersects the two diagonals orthogonally at their midpoints (\Figref{fig:Bennett_angles}, see, e.g., \cite[p.~555]{Quadrics}).
The same holds for spherical isograms when the sides are replaced by chords.
Consequently, their axis $m$ is the diameter of $\Sphere$ that connects the two intersections $S$ and $S^\ast$ of the diagonal great circles.
For non-crossed isograms the axis meets the diagonal arcs at their midpoints.
Otherwise, $S$ and $S^\ast$ are spherical centers of the great circle $s$ that meets the diagonal arcs orthogonally at their midpoints (note \Figref{fig:sph_isogram}, left and right).
Due to these symmetries, opposite angles at spherical isograms are congruent.

If one vertex of a spherical isogram is replaced by its antipode, then we obtain a quadrangle where the lengths of opposite sides sum up to $\pi$. 

\begin{lemma}\label{lem:symm_center}
For any two different oriented great circles $g_1, g_2$ on $\Sphere$ there exists a unique pair of antipodal centers of symmetry $S, S^\ast \in\Sphere$, i.e., spherical centers of a halfturn that maps $g_1$ to $g_2$, and vice versa.
The centers are located on the common perpendicular of $g_1$ and $g_2$.
Also the reflection in the great circle $s$ with spherical centers $S$ and $S^\ast$ exchanges the oriented great circles $g_1$ and $g_2$. 
\end{lemma}

\begin{Proof}
The planes $\gamma_1, \gamma_2$ spanned by the great circles $g_1, g_2$ intersect each other along a line $d$ through $O$.
This line $d$ must be preserved by the wanted halfturn. 
Therefore, its axis lies in the diameter plane orthogonal to $d$ and moreover in one of the two planes that bisect the angle between $\gamma_1$ and $\gamma_2$.
\\
The images of $g_1$ under the two possible halfturns have different orientation.
Hence, just one of them yields the wanted orientation of $g_2$.
The same holds for the reflection of $g_1$ in the diameter plane that is orthogonal to the axis of the correct halfturn. 
\end{Proof}

From now on we assume that the four sides of a spherical isogram are basis, coupler and arms of a spherical coupler motion (\Figref{fig:sph_isogram}).
It is well known that this special case of a coupler motion splits into two rational one-parameter motions with bifurcations in the aligned poses.
During each rational motion the arms rotate in such a way that the tangents of the halved rotation angles $t_i:= \tan \frac{\phi_i}{2}$ are proportional (see, e.g., \cite{Bricard} or \cite[eq.~(9)]{Sta_128}).
Referring to the notions used in \Figref{fig:Bennett_angles} holds
\begin{equation}\label{eq:proportional}
   t_2 = c_{21}t_1, \ \mbox{where} \ c_{21} = \frac{\sin(\alpha - \beta)}{\sin\alpha \pm \sin\beta}, \ \alpha = \wh{R_{01}R_{02}}, \ \beta = \wh{R_{01}R_{31}}.
\end{equation}   
 
\begin{figure}[htb] 
  \centerline{\includegraphics[width=70mm]{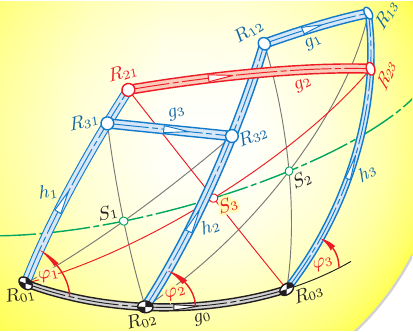}}
  \caption{This triple of spherical isograms is flexible and has aligned centers of symmetry (\Lemref{lem:compound_isograms}).}
  \label{fig:triple}
\end{figure}

\begin{lemma}\label{lem:compound_isograms} 
On a sphere $\Sphere$, let $R_{01}R_{02}$ and $R_{02}R_{03}$ be the aligned bases of two isograms where the arms through $R_{02}$ are aligned, as well (\Figref{fig:triple}).
Then, the simultaneous continuous movement of these two isograms is compatible with a third isogram with the basis $R_{01}R_{03}$ and with arms respectively aligned with arms of the given isograms.
\\[0.5mm]
The centers of symmetry $S_1$, $S_2$ and $S_3$ of the three isograms are aligned.
\end{lemma}

\begin{Proof}
Let $\phi_1$, $\phi_2$ and $\phi_3$ denote the angles of rotation of the three arms relative to the initial aligned pose on the great circle $g_0$ (\Figref{fig:triple}).
By \eqref{eq:proportional} the respective tangents $t_1$ and $t_2$ of halved rotation angles as well as $t_2$ and $t_3$ are proportional. 
Thus, the composition results again in an isogram motion. 
In order to determine the coupler of the resulting isogram, one can use the formulas for the factors of proportionality.

\noindent
We prefer a more geometric solution:
Given any non-aligned pose, let $h_i$ be the great circle spanned by the arm through  $R_{0i}$ for $i=1,2,3$. 
Additionally, we transfer an arbitrarily chosen orientation of the great circle $g_0$ by the given rotations to an orientation of $h_1$, $h_2$ and $h_3$.
\\
Referring to \Lemref{lem:symm_center}, we look for the center $S_3$ of a halfturn $\sigma_3$ that sends $h_1$ to the opposite orientation $\ol h_3$ of $h_3$. 
This center $S_3$ is located on the great circle that meets $h_1$ and $h_3$ orthogonally.
The halfturn about $S_3$ sends the basis $R_{01}R_{03}$ to the wanted coupler $R_{23}R_{21}$ and the great circle $g_0$ to the opposite orientation of the great circle $g_2$ that extends the coupler.
\\
On the other hand, there exist centers $S_1$ and $S_2$ of halfturns $\sigma_1$ and $\sigma_2$ that take the bases $R_{01}R_{02}$ and $R_{02}R_{03}$ into the respective couplers $R_{32}R_{31}\subset g_3$ and $R_{13}R_{12}\subset g_1$.
We summarize:
\setlength{\arraycolsep}{1.0mm}
\begin{equation}\label{eq:sigma_1..3}
 \begin{array}{lllll}
  \sigma_1\!: &g_0 \leftrightarrow \ol g_3, &h_1\leftrightarrow \ol h_2,
   &R_{01}\leftrightarrow R_{32}, &R_{02}\leftrightarrow R_{31},
  \\[0.4mm] 
  \sigma_2\!: &g_0 \leftrightarrow \ol g_1, &h_2\leftrightarrow \ol h_3,
   &R_{02}\leftrightarrow R_{13}, &R_{03}\leftrightarrow R_{12},
  \\[0.4mm] 
  \sigma_3\!: &g_0 \leftrightarrow \ol g_2, &h_3\leftrightarrow \ol h_1,
   &R_{03}\leftrightarrow R_{21}, &R_{01}\leftrightarrow R_{23}.
 \end{array}
\end{equation}

\noindent
The composition of $\sigma_1$ with $\sigma_2$, which is denoted by $\sigma_2\sigma_1$, is a rotation $\rho_{21}$ about the diameter $p$ orthogonal to the axes of $\sigma_1$ and $\sigma_2$ through twice of the angle $\wh{S_1S_2}$ and sends $h_1$ to $h_3$. 
The reversed product $\sigma_1\sigma_2$ gives the inverse displacement $\rho_{21}^{-1}$, i.e., the rotation about $p$ through the same angle in the opposite direction.

The composition $\sigma_3 \sigma_2 \sigma_1$ of the three halfturns is an orientation preserving displacement, consequently a rotation that reverses the orientation of the great circle $h_1$.
This must be a halfturn about a diameter of the great circle $h_1$, and this is involutive.
Therefore holds
\[  \sigma_3 \sigma_2 \sigma_1 = \left(\sigma_3 \sigma_2 \sigma_1\right)^{-1} 
   = \sigma_1 \sigma_2 \sigma_3. 
\]
The equivalent equation
\begin{equation}\label{eq:transform_rotation}
   \sigma_3\klzwi (\sigma_2 \sigma_1)\klzwi \sigma_3 = \sigma_1 \sigma_2
\end{equation}
reveals that the halfturn $\sigma_3$ transforms the rotation $\rho_{21}$ into $\rho_{21}^{-1}$, which means, that it reverses the axis $p$ of $\rho_{21}$.
Hence, the axis of $\sigma_3$ meets $p$ orthogonally at the sphere's center $O$ and is therefore coplanar with the axes of $\sigma_1$ and $\sigma_2$.  
\end{Proof}

According to \eqref{eq:sigma_1..3}, the points $R_{01}$ and $R_{32} = \sigma_1(R_{01})$ have equal distances to the plane spanned by the axes of $\sigma_1$, $\sigma_2$ and $\sigma_3$.
And there hold similar relations for the other vertices.


\begin{figure}[htb] 
  \centerline{\includegraphics[width=110mm]{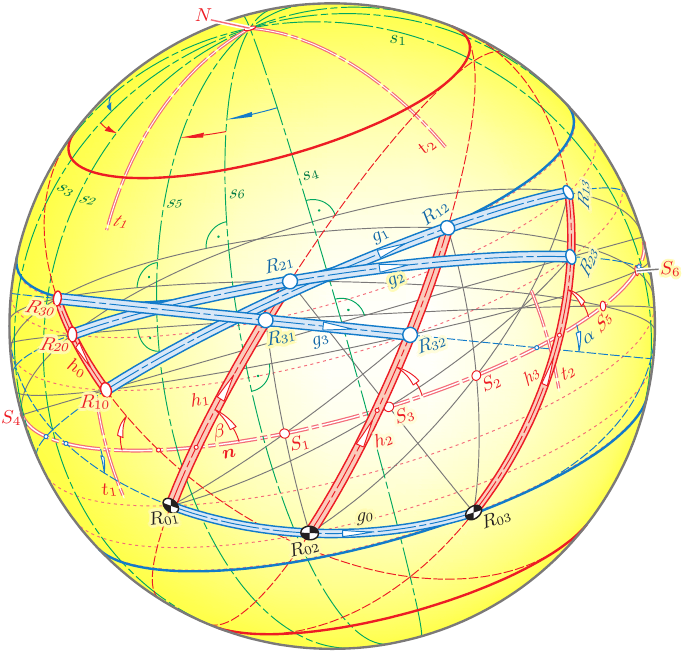}}
  \caption{Bennett's spherical 8-bar linkage is the union of six isograms with centers of symmetry $S_1,\dots,S_6$ aligned on the great circle $n$ with spherical center $N$ (\Lemref{lem:compound_isograms}).
The 12 joints $R_{ij}$ are located by four on pairs of circles parallel and symmetric with respect to $n$.
The 8 bars can be subdivided into two quadruples that make respectively equal angles with $n$.}
  \label{fig:sph_linkage2}
\end{figure}

\smallskip
Now we are able to extend the triple of isograms depicted in \Figref{fig:triple} by additional triples to obtain Bennett's spherical 8-bar linkage shown in \Figref{fig:sph_linkage2}:

The great circle $h_3$ carries the bases $R_{03}R_{13}$ and $R_{03}R_{23}$ of two isograms with the centers $S_2$ and $S_3$. 
By virtue of \Lemref{lem:compound_isograms}, there exists another spherical isogram with the basis $R_{13}R_{23}$, with arms along $g_1$ and $g_2$, and with the center of symmetry $S_4$ that must be aligned with $S_2$ and $S_3$.

Similarly, the bases $R_{01}R_{31}$ and $R_{01}R_{21}$ on $h_1$ belong to the isograms with centers $S_1$ and $S_3$ and give rise to another isogram with basis $R_{21}R_{31}$, with arms along $g_2$ and $g_3$, and with the center $S_5$ aligned with $S_1$ and $S_3$.
Finally, the bases $R_{02}R_{12}$ and $R_{02}R_{32}$ on $h_2$ define a sixth isogram.  It has the basis $R_{12}R_{32}$, the arms along $g_1$ and $g_3$, and the center $S_6$ aligned with all other centers $S_1,\dots,S_5$.
We summarize so far:
\begin{equation}\label{eq:sigma_4..6}
    \sigma_4\!: \,g_1 \leftrightarrow \ol g_2, \quad
    \sigma_5\!: \,g_2 \leftrightarrow \ol g_3, \quad
    \sigma_6\!: \,g_3 \leftrightarrow \ol g_1.
\end{equation}
It remains to prove that the couplers of the three additional isograms lie on the same great circle $h_0$ (note \Figref{fig:sph_linkage2}):

The great circles spanned by the first two additional couplers arise, respectively, from $h_3$ by the halfturn $\sigma_4$ about $S_4$ and from $h_1$ by the halfturn $\sigma_5$ about $S_5$.
By \eqref{eq:sigma_1..3}, $h_3$ and $h_1$ are images of the opposite orientation $\ol h_2$ of $h_2$ under $\sigma_2$ and $\sigma_1$, respectively.
This implies that the couplers of the first two additional isograms are located on the great circles $\sigma_4\sigma_2(\klzwi\ol h_2)$ and on $\sigma_5\sigma_1(\klzwi\ol h_2)$.  
These two great circles are equal since we can state 
\begin{equation}\label{eq:42=51}
  \rho_{42}:= \sigma_4\sigma_2 = \sigma_5\sigma_1 =:\rho_{51}. 
\end{equation}

\begin{Proof}
Both products are rotations about the axis $p$ orthogonal to the plane of the great circle through the six centers of symmetry.
Both products send $g_0$ to the right orientation of $g_2$, the first via $\ol g_1$, the latter via $\ol g_3$.     
Consequently, if we apply to $g_0$ the first rotation and afterwards the inverse of the second rotation, we obtain again $g_0$ with the correct orientation.
Since $p$ cannot be the axis of the great circle $g_0$, the composition of the two rotations must be the identity, which confirms the equation \eqref{eq:42=51}.
\end{Proof}

In a similar way we can prove that the sixth coupler, which is placed on $\sigma_6(h_2)$, is aligned with that of the fifth isogram placed on $\sigma_5(h_1)$.
We obtain it as a consequence of $h_2 = \sigma_2(\ol h_3)$, $h_1 = \sigma_3(\ol h_3)$ along with
\begin{equation}\label{eq:62=53}
  \rho_{62}:= \sigma_6\sigma_2 = \sigma_5\sigma_3 =:\rho_{53} \quad \mbox{and} \quad
  \rho_{61}:= \sigma_6\sigma_1 = \sigma_4\sigma_3 = \rho_{43}. 
\end{equation} 
The products in \eqref{eq:42=51} and \eqref{eq:62=53} are rotations about the axis $p$ through the point $N\in\Sphere$, and we conclude
\begin{equation}\label{eq:rotations}
  \rho_{61}\!: g_0\mapsto g_1,\ h_1\mapsto h_0,\quad
  \rho_{42}\!: g_0\mapsto g_2,\ h_2\mapsto h_0,\quad
  \rho_{53}\!: g_0\mapsto g_3,\ h_3\mapsto h_0.
\end{equation}  
They imply further identities like
\begin{equation}\label{eq:rotations2}
  \rho_{54}:= \sigma_1\sigma_2:= \rho_{12}, \quad
  \rho_{65}:= \sigma_2\sigma_3:= \rho_{23}, \quad
  \rho_{46}:= \sigma_3\sigma_1:= \rho_{31},  
\end{equation}
which result in congruent angles at $N$ between the symmetry axes $s_i$ of isograms as marked in \Figref{fig:sph_linkage2}.   
In other words, the pairs $(s_1,s_4)$, $(s_2,s_5)$ and $(s_3,s_6)$ share the spherical angle bisectors $t_1$ and $t_2$. 
By virtue of \eqref{eq:rotations}, the distances of $R_{10}, R_{01},R_{23}$, and $R_{32}$ to the plane spanned by the centers $S_1\dots,S_6$ are equal.
This confirms the local symmetries that were already taken in consideration by Bennett in \cite[Figs.~5 and 6]{Bennett}.

\begin{theorem}\label{thm:1} 
Bennett's spherical six-bar linkage as depicted in \Figref{fig:spher_linkage1} is continuously flexible with degree of freedom $\mbox{d.o.f.} = 1$.
In each pose, the centers of symmetry $S_1,\dots,S_6$ of the contained isograms are placed on a great circle $n$.
\\[0.5mm]
For $i=1,2,3$ the rotation with $g_0\mapsto g_i$ about the spherical center $N$ of $n$ takes $h_i$ to $h_0$. 
The intersections of $n$ with the extended sides $g_0,\dots,g_3$ and $h_0,\dots,h_3$ are respectively symmetric w.r.t.\ orthogonal great circles $t_1,t_2$ through $N$ (\Figref{fig:sph_linkage2}). 
\end{theorem}

\begin{Proof}
As mentioned above, the product $\tau_{321}:= \sigma_3\sigma_2\sigma_1$ is a halfturn that reverses the orientation of $h_1$.
Therefore, the axis of $\tau_{321}$ is a diameter of $h_1$ and, on the other hand, coplanar with the axes of $\sigma_1,\dots,\sigma_6$.
Reflection in $t_1$ or $t_2$ reveals that $\tau_{654}:= \sigma_6\sigma_5\sigma_4$
is a halfturn about an axis symmetric to that of $\tau_{321}$, while on the other hand $\tau_{654}$ changes the orientation of $g_1$ by virtue of \eqref{eq:sigma_4..6}. 
After similar conclusions for other threefold products, this confirms the last statement.
\end{Proof}

\section{The dualized version}
\begin{figure}[hbt] 
  \centerline{\includegraphics[width=90mm]{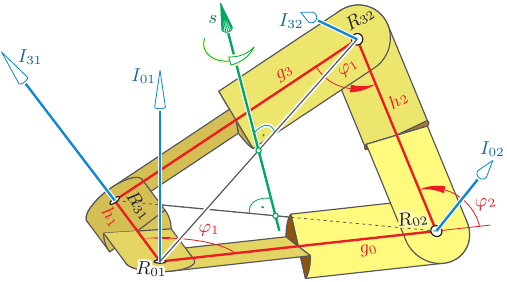}}
  \caption{The Bennett isogram is a flexible 4R loop where the common perpendiculars between adjacent hinges form an isogram.
The reflection in the line $s$ exchanges opposite sides, angles and hinges.}
  \label{FIG_bennet2} 
\end{figure}

The spatial analogues of spherical isograms are Bennett isograms, i.e., flexible 4R chains where the common perpendiculars of neighboring hinges form a closed skew isogram (\Figref{FIG_bennet2}).
In each spatial pose, there exists a line $s$ such that the reflection in $s$ exchanges opposite vertices, hinges and angles (\Figref{fig:Bennett_angles}). 
According to \cite{Brunn}, the Bennett isogram is the only 4R loop which is continuously flexible. 

\begin{figure}[htb] 
  \centerline{\includegraphics[width=90mm]{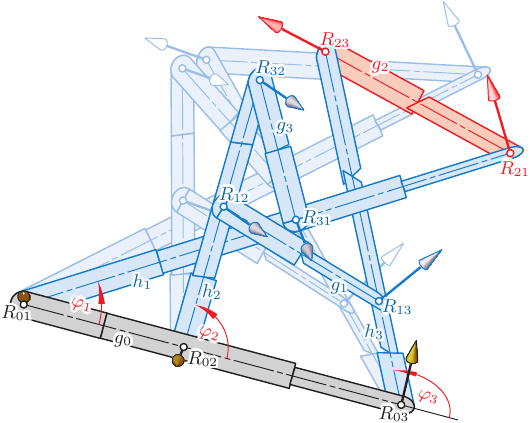}}
  \caption{A flexible triple of Bennett isograms with aligned bases on $g_0$, arms on the lines $h_1$, $h_2$ and $h_3$, and couplers on the lines $g_3$, $g_1$ and $g_2$.}
  \label{fig:spatial_linkage}
\end{figure}

The hinges of a Bennett isogram define as their spherical image a spherical isogram.
This reveals that the angles of rotation for neighbouring sides obey again the equation \eqref{eq:proportional}.
Hence, similar to \Lemref{lem:compound_isograms} the composition of two Bennett isograms with aligned bases $R_{01}R_{02}$ and $R_{01}R_{02}$ on $g_0$ and aligned arms $R_{02}R_{32}$ and $R_{02}R_{32}$ on $h_2$, as shown in \Figref{fig:spatial_linkage},  yield a transmission according to \eqref{eq:proportional} from the first arm $R_{01}R_{31}$ on $h_1$ of the first isogram to the last arm $R_{01}R_{13}$ on $h_3$ of the second.
This shows that a third Bennett isogram with basis $R_{01}R_{03}$ is compatible with the motions of the first two isograms.

When we assign orientations to the sides $g_0,h_1,h_2,h_3$ like previously, then the reflection in the axis of symmetry $s_3$ of the third isogram is defined by sending $h_1$ to the reversed orientation $\ol h_3$ of $h_3$. 
Therefore, it passes through the midpoint of the shortest segment between $h_1$ and $h_3$ and has the direction of an angle bisector between the directions of $h_1$ and $\ol h_3$.

We adopt the notation from the spherical case and denote the line reflections of the three isograms with $\sigma_1$, $\sigma_2$ and $\sigma_3$.
Thus, we can follow the lines of the foregoing section nearly word for word.
The only difference is that now the composition $\sigma_2\sigma_1$ of the first two reflections is a helical displacement about the common perpendicular $n$ of the two respective axes $s_1$ and $s_2$, that takes $h_1$ with $R_{31}$ and the hinge $I_{31}$  to $h_3$ with $R_{13}$ and $I_{13}$.
Since the product $\tau_{321}:= \sigma_3\sigma_2\sigma_1$ is a helical displacement that keeps the line $h_1$ fixed but reverses its orientation, $\tau_{321}$ must be the reflection in a line orthogonal to $h_1$.
This implies by \eqref{eq:transform_rotation} like above, that $s_3$ intersects $n$ perpendicularily, too.
Moreover, the axis of $\tau_{321}$ is the common perpendicular of $h_1$ and $n$.  

\begin{figure}[htb] 
  \centerline{\includegraphics[width=105mm]{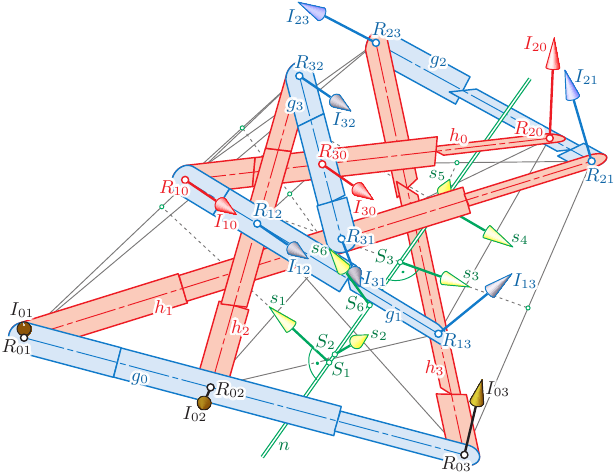}}
  \caption{The spatial version of Bennett's spherical 8-bar linkage consists of 8 bodies $\Gamma_0,\dots,\Delta_3$ and 12 hinges $I_{01}\dots I_{32}$, arranged in 6 Bennett isograms.
In each spatial pose, the axes of symmetry $s_1,\dots,s_6$ of the six isograms have the line $n$ as their common perpendicular.}
  \label{fig:spatial_linkage_kmplt}
\end{figure}

\begin{figure}[htb] 
  \centerline{\includegraphics[width=105mm]{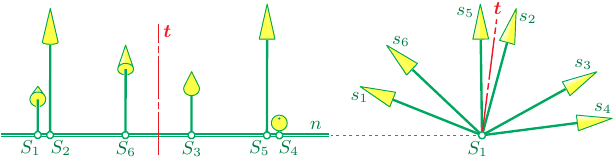}}
  \caption{The axes of symmetry $s_1,\dots,s_6$ meet the line $n$ orthogonally.
  They are placed in such a way that the pairs $(s_1,s_4)$, $(s_2,s_5)$ and $(s_3,s_6)$ have a common axis $t$ of symmetry.   
Left: View orthogonal to $n$. Right: Point view of $n$.}
  \label{fig:perp_axes} 
\end{figure}

After extending the first triple of Bennett isograms by three additional isograms, we obtain additional line reflections $\sigma_4, \sigma_5$ and $\sigma_6$ that satisfy the equations \eqref{eq:42=51}---\eqref{eq:rotations2}.
This confirms that the additional isograms have aligned couplers on the line $h_0$ while the axes of symmetry $s_1\dots s_3$ of all six isograms are placed such that the pairs $(s_1,s_4)$, $(s_2,s_5)$ and $(s_3,s_6)$ have a common axis of symmetry $t$ orthogonal to $n$ (\Figref{fig:perp_axes}).
As a consequence, the lines $g_0,\dots,g_3$ have equal distances to $n$ and make congruent angles with $n$, and the same holds for the distances and angles between the lines $h_0,\dots,h_3$ and $n$.

\begin{theorem}\label{thm:2} 
The spatial version of Bennett's spherical eight-bar linkage is a compound of six Bennett isograms (\Figref{fig:spatial_linkage_kmplt}) and continuously flexible with $\mbox{d.o.f.} = 1$.
In each pose, the axes of symmetry $s_1,\dots,s_6$ of the contained isograms have a common perpendicular $n$.
\\[0.5mm]
For $i=1,2,3$ holds that the helical displacement about $n$ with $g_0\mapsto g_i$ takes simultaneously $h_i$ to $h_0$. 
The common perpendiculars between $n$ and $g_0,\dots,g_3$ are respectively symmetric to the common perpendiculars between $n$ and $h_0,\dots,h_3$ with respect to an axis $t$ (\Figref{fig:perp_axes}).   
\end{theorem}

\paragraph{Notes and Comments}\hfill{}\newline
The author expresses sincere thanks to Hans-Peter Schr\"ocker, University of Innsbruck,  for useful advice.

%
%

\end{document}